\RequirePackage{fix-cm}

\documentclass[twocolumn]{svjour3}          % twocolumn
\smartqed  % flush right qed marks, e.g. at end of proof
\usepackage{booktabs}
\usepackage{subcaption}
\usepackage{graphicx}
\usepackage{epstopdf}
\usepackage{amsmath,amsfonts}
\usepackage{commath}
\usepackage{todonotes}

\usepackage{tikz}
\usepackage{pgfplots}
\pgfplotsset{compat=1.13}
\usetikzlibrary{plotmarks,calc,math,shapes}

\usepackage[utf8]{inputenc}

\renewcommand{\Re}{\mathit{Re}}
\newcommand{\pu}{{\text{p.u.}}}
\newcommand{\va}{\phi} % variable for voltage angle
\newcommand{\fuel}{\varepsilon}

\newcommand\RR{{\mathbb{R}}}
\newcommand\Nu{{\mathit{Nu}}} % number of control variables
\newcommand\Ny{{\mathit{Ny}}} % number of state variables
\newcommand\td[1]{\frac{\mathit{d}}{\mathit{d} #1}}
\renewcommand\pd[1]{\frac{\partial}{\partial #1}}

\begin{document}

\title{Optimal control of compressor stations in a coupled gas-to-power network\thanks{The authors gratefully thank the BMBF project ENets (05M18VMA) for the financial support.}
%Grants or other notes
}

%\titlerunning{Coupled gas-to-power network}        % if too long for running head

\author{Eike Fokken \and Simone Göttlich \and Oliver Kolb}


\institute{Eike Fokken \at
              University of Mannheim \\
              \email{fokken@uni-mannheim.de}           %  \\
           \and
           Simone Göttlich \at
           University of Mannheim \\
           \email{goettlich@uni-mannheim.de}           %  \\
           \and
           Oliver Kolb \at
           University of Mannheim \\
           \email{kolb@uni-mannheim.de}           %  \\
}


\maketitle

\begin{abstract}
We introduce a tool for simulation and optimization of gas pipeline networks coupled to power grids by gas-to-power plants. 
The model under consideration consists of the isentropic Euler equations to describe the gas flow coupled
to the AC powerflow equations. A compressor station is installed to control the gas pressure such that certain bounds are satisfied.
A numerical case study is presented that showcases effects of fast changes in power demand on gas pipelines and necessary operator actions.
  %35L65, 65M08
\end{abstract}

\section{Introduction}
\label{intro}
Renewable power sources have an ever increasing share of all power sources.
Though renewable energy has been developed in recent years with great success, its intermittent and unpredictable nature raises 
the difficulty to balance the energy production and consumption \cite{Gahleitner2013,ZengFangLiChen2016}.
A frequent proposal is to use gas turbine plants to compensate for sudden drops in power of renewable sources
because these plants are relatively flexible in comparison to coal or nuclear plants.
A welcome advantage of this approach is the possibility to run gas plants with fuel produced from renewable electricity via power-to-gas plants, thereby reducing or even negating carbon emissions of the gas plants.  For a review of power-to-gas capability see \cite{Gahleitner2013}.

It is desirable to have a joint optimal control framework for power and gas sector of the energy system to model this compensation.
So far only steady-state flow in the gas network has been considered \cite{CHERTKOV2015541,ZengFangLiChen2016,Zlotnik2016},
which may be too coarse for several applications.
Therefore, we focus on an optimal control strategy for the instationary gas network model~\cite{banda_herty_klar} coupled to a power grid~\cite{Bienstock} via compressor stations,
see for example \cite{Herty2007,Mak2018,CompressorModels2017}. 
The mathematical foundation for the gas-to-power coupling has been recently introduced in \cite{fokken_goettlich_kolb}, 
where conditions for the well-posedness have been derived and proved.
This work is the first attempt to model this interaction and yields an understanding of the underlying equations.  The next step will be real-world scenarios.

\section{Optimal control problem}
\label{sec:model}

The gas dynamics within each pipeline of the considered gas networks are modeled by the isentropic Euler equations, supplemented with suitable coupling and boundary conditions. For the power grid, we apply the well-known powerflow equations. The coupling between gas and power networks at gas-driven power plants is modeled by (algebraic) demand-dependent gas consumption terms.
To react on the demand-dependent influences on the gas network, controllable devices as compressor stations are considered within the gas network. 
The aim is to fulfill given state restrictions like pressure bounds whereas at the same time the entire fuel gas or power consumption of the compressor stations is to be minimized. 

\begin{figure}[h!]
  \centering
 \includegraphics[width=0.4\textwidth]{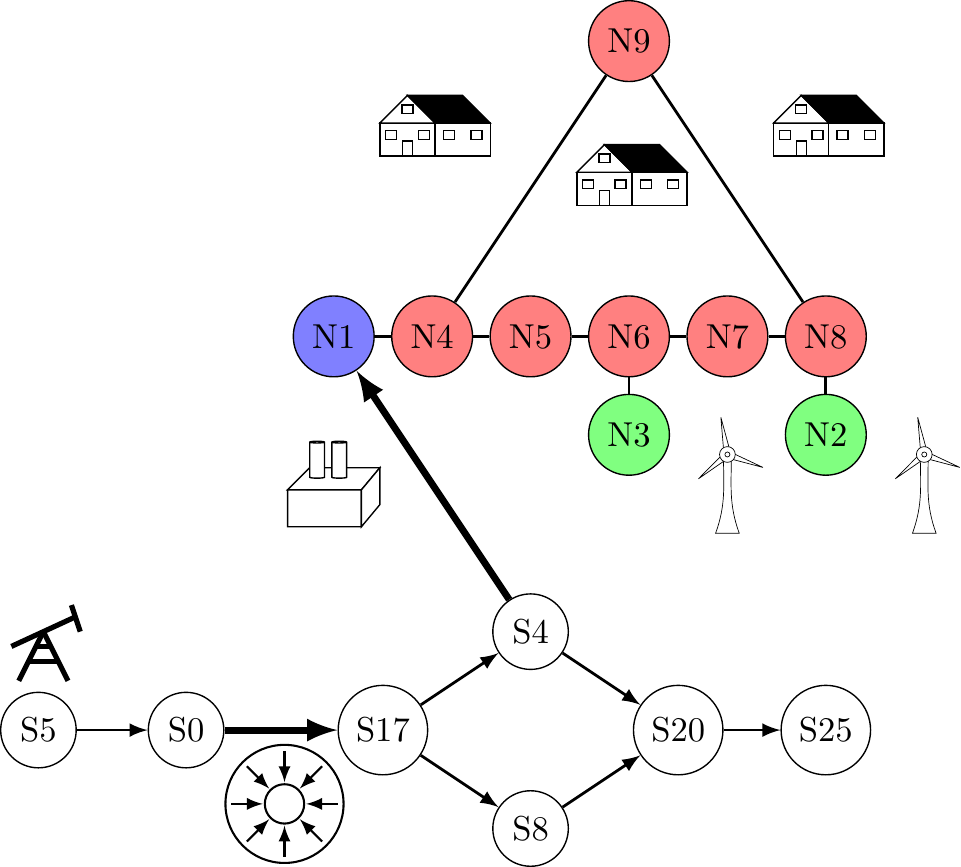}
 \caption{Gas network connected to a power grid. Red nodes are \textcolor{red!50}{PQ/demand nodes}, green nodes are generators (\textcolor{green!50}{PV nodes}) and the blue node is the \textcolor{blue!50}{slack bus} (also a generator, with gas consumption of the form $\fuel(P)=a_0 + a_1 P + a_2 P^2$). The circle symbol indicates a compressor station.}
 \label{fig:gas2power}
\end{figure}

The network under consideration is similar to the one presenetd in \cite{fokken_goettlich_kolb} and is depicted in Figure \ref{fig:gas2power}.
In addition, there is now a compressor with a time-dependent control $u(t)$, which is used to satisfy pressure bounds.
The optimal control problem is to minimize compressor costs while satisfying power demand and gas dynamics:
\begin{equation*} 
 \begin{array}{cl}
  \min_{u(t)} & \text{compressor costs}\\
  \text{s.t.} & \text{isentropic Euler equations (Subsection \ref{sec:gas-model})}\\ 
	& \text{gas coupling conditions (Subsection \ref{sec:coupling-at-gas})}\\ 
	& \text{compressor equation (Subsection \ref{sec:compressor-stations})}\\
  & \text{powerflow equations (Subsection \ref{sec:power-model})}\\ 
  & \text{gas-power-coupling (Subsection \ref{sec:coupling})}\\
  & \text{pressure bounds.}
 \end{array}
\end{equation*}

Mathematically, this is an instationary nonlinear optimization problem constrained by partial differential equations, see \cite{Hinze2009} for an overview.
To solve the problem, we make use of a first-discretize-then-optimize approach and apply the interior point solver IPOPT \cite{Ipopt}. The necessary gradient information for IPOPT, i.e., gradients with respect to all controllable devices, is efficiently computed via adjoint equations. Here, the underlying systems can be solved time-step-wise (backwards in time), where additionally the sparsity structure is exploited.
  
We remark that the cost function is given in Subsection \ref{sec:compressor-stations}, while
the bounds on the pressure are introduced as box constraints within the numerical optimization procedure in Section \ref{sec:2}.
Within the Subsections \ref{sec:gas-model}--\ref{sec:coupling}, 
we now describe the constraints of the optimal control problem in detail and focus on the technical details.
\subsection{The isentropic Euler equations}
\label{sec:gas-model}
The gas network is modeled by the isentropic Euler equations (see \cite{Colombo2008,banda_herty_klar}), which govern gas flow in each pipeline between nodes,
  \begin{equation}\label{eq:4}
    \begin{pmatrix}\rho\\  q\end{pmatrix}_t + \begin{pmatrix} q\\ p(\rho) + \frac{ q^2}{\rho}\end{pmatrix}_x = \begin{pmatrix}0\\S(\rho, q) \end{pmatrix},
\end{equation}
where $\rho$ is the density, $q= \rho v$ is the flow, $p(\rho) = \kappa\rho^\gamma$ is the pressure function.  In our example we use $\gamma =1$ and $\kappa = c^2 = (340\frac{\text{m}}{\text{s}})^2$, that is, the isothermal Euler equations with speed of sound $c$.  The Euler equations must be met for $0\leq x \leq l_e$ and $t \geq 0$, where $l_e$ is the length of the pipe.
Furthermore, $S$ is a source term,
\[
 S(\rho,q) = -\frac{\lambda(q)}{2d_e} \frac{q \vert q\vert}{\rho},
\]
where $\lambda$ is the solution to the Prandtl-Colebrook formula,
\begin{equation*}
 \frac{1}{\sqrt{\lambda}} = -2\log_{10} \left( \frac{2.51}{\Re(q) \sqrt{\lambda}} + \frac{k_e}{3.71 d_e} \right).
\end{equation*}
 The Reynolds number $\Re$ is given by
 \begin{equation*}
 \Re(q) = \frac{d_e}{\eta} q
\end{equation*}
with dynamic viscosity
\begin{equation*}
 \eta_{\phantom{e}} = 10^{-5} \frac{\text{kg}}{\text{ms}}.
\end{equation*}
The roughness $k_e$ and diameter $d_e$ of the pipes are all the same in our example and given by
\begin{equation*}
  \begin{aligned}
    k_e &= 5\cdot 10^{-4}\,\text{m},\\
    d_e &= 6\cdot 10^{-1}\,\text{m}.
\end{aligned}
\end{equation*}

\subsection{Coupling at gas nodes}
\label{sec:coupling-at-gas}
At the nodes we use the usual Kirchhoff-type coupling conditions: The pressure is the same near the node in all pipes connected to it and the flows must add up to zero (where the sign for inflow is positive, the sign for outflow negative),
\begin{subequations}\label{eq:couplings}
  \begin{align}
    p_\text{pipe} &= p_\text{node}\label{eq:pressure_coupling},\\
    \sum_\text{ingoing pipes} q_\text{pipe} &= \sum_\text{outgoing pipes} q_\text{pipe}.\label{eq:flow_coupling}
  \end{align}
\end{subequations}

The example gas network we are using is a small part of the GasLib-40 network~\cite{Humpola_et_al:2015}.  In Table~\ref{tab:gasNetwork} the only remaining parameter, the length $l_e$ of each pipe, is gathered. Note that no length is given for the arc connecting S$0$ and S$17$, because only a compressor is situated between these nodes.
\begin{table}[hbt]
  \centering
  \caption{Parameters of the gas network\label{tab:gasNetwork}}
  \begin{tabular}{@{}lllc}
    \toprule
    Pipe & From & To & Length $l_e$ [km] \\
    \midrule
    P10 & S4  & S20 & 20.322  \\
    P20 & S5  & S0  & 20.635  \\
    P21 & S17 & S4  & 10.586  \\
    P22 & S17 & S8  & 10.452  \\
    P24 & S8  & S20 & 19.303  \\
    P25 & S20 & S25 & 66.037  \\
    \bottomrule
  \end{tabular}
\end{table}
\subsection{Compressor stations}
\label{sec:compressor-stations}
To compensate for pressure losses in the gas network, we consider compressor stations, which are also modelled as arcs.
Those arcs have (time-depen\-dent) in- and outgoing pressure ($p_\text{in}$, $p_\text{out}$) and flux values ($q_\text{in}$, $q_\text{out}$).
In general, the two separate flux values allow the modelling of fuel gas consumption of the compressor station, whereas we will consider an external power supply for the compressor and therefore have $q_\text{in}=q_\text{out}$.
The power consumption is modelled as a quadratic function of the power required for the compression process.
Denoting this as function $P(p_{\text{in}}(t),q_{\text{in}}(t),p_{\text{out}}(t),q_{\text{out}}(t))$, our objective function in the optimal control problem is of the form
\begin{equation}\label{eq:compressorCosts}
 \int\limits_0^T P(p_{\text{in}}(t),q_{\text{in}}(t),p_{\text{out}}(t),q_{\text{out}}(t)) dt.
\end{equation}
Note that the power consumption does not influence the network dynamics and is therefore only of interest for the optimization procedure.
For our investigations below it is sufficient to know that the consumption and therewith the costs increase if the ratio $p_\text{out}/p_\text{in}$ increases.
For the details of the power consumption model, we refer to~\cite[Section~3.2.3]{CompressorModels2017}.
The influence of the compressor station on the network dynamics is modelled by the control $u(t)$ of the pressure difference:
\begin{equation*}
 p_\text{out}(t) - p_\text{in}(t) = u(t).
\end{equation*}

\subsection{Power model}
\label{sec:power-model}
For the power grid we use the AC powerflow equations (see \cite{grainger2016power} for an introduction),
\begin{equation}\label{eq:3}
  \begin{aligned}
    P_{k}  & = \sum_{j=1}^N \abs{V_k}\abs{V_j}(G_{kj}\cos(\phi_{k,j})+B_{kj}\sin(\phi_{k,j})),\\
    Q_{k}  & = \sum_{k=1}^N   \abs{V_k}\abs{V_j}(G_{kj}\sin(\phi_{k,j})-B_{kj}\cos(\phi_{k,j})),
  \end{aligned}
\end{equation}
where $P_k$, $Q_k$ are real and reactive power at node $k$, $\abs{V_k}$ is the voltage amplitude, $\phi_k$ is the phase (and $\phi_{k,j} = \phi_k-\phi_j$) and $B_{kj}$, $G_{kj}$ are parameters of the transmission lines between nodes $k$ and $j$ or of the node $k$ for $B_{kk}$ and $G_{kk}$.
 Each node is either the slack bus (in our case N$1$; $V$ and $\phi$ given), a generator bus (N$2$ and N$3$; $V$ and $P$ given) or a load bus (N$4$ through N$9$; $P$ and $Q$ given).
 All in all for $N$ nodes we get $2N$ equations for $2N$ variables.

The considered power grid is taken from the example ``case9'' of the MATPOWER Matlab programming suite~\cite{MATPOWER}.  A per-unit system is used, whose base power and voltage are 100MW and 345kV respectively.  The corresponding node and transmission line parameters are gathered in Table~\ref{tab:powerGrid}, these are the entries of the nodal admittance matrix (see \cite{grainger2016power}).
\begin{table}[hbt]
  \centering
   \caption{Parameters of the power grid (p.u.)\label{tab:powerGrid}}

   \begin{subtable}{0.5\textwidth}
    \centering
    \subcaption{Busses}
    \begin{tabular}{@{}lll}
   \toprule
  Node & $G$ & $B$ \\
  \midrule
  N1 & 0.0000 & -17.3611 \\
  N2 & 0.0000 & -16.0000 \\
  N3 & 0.0000 & -17.0648 \\
  N4 & 3.3074 & -39.3089 \\
  N5 & 3.2242 & -15.8409 \\
  N6 & 2.4371 & -32.1539 \\
  N7 & 2.7722 & -23.3032 \\
  N8 & 2.8047 & -35.4456 \\
  N9 & 2.5528 & -17.3382 \\
  \bottomrule
    \end{tabular}
 \end{subtable}
 
 \begin{subtable}{0.5\textwidth}
   \centering
     \subcaption{Transmission lines}
     \begin{tabular}{@{}lllll}
    \toprule
  Edge & From & To & $G$ & $B$ \\
  \midrule
  TL14 & N1 & N4 & \phantom{-}0.0000 & 17.3611 \\
  TL45 & N4 & N5 & -1.9422 & 10.5107 \\
  TL56 & N5 & N6 & -1.2820 & \phantom{1}5.5882 \\
  TL36 & N3 & N6 & \phantom{-}0.0000 & 17.0648 \\
  TL67 & N6 & N7 & -1.1551 & \phantom{1}9.7843 \\
  TL78 & N7 & N8 & -1.6171 & 13.6980 \\
  TL82 & N8 & N2 & \phantom{-}0.0000 & 16.0000 \\
  TL89 & N8 & N9 & -1.1876 & \phantom{1}5.9751 \\
  TL94 & N9 & N4 & -1.3652 & 11.6041 \\
  \bottomrule
     \end{tabular}
   \end{subtable}
 \end{table}

\subsection{Coupling}
\label{sec:coupling}
The last ingredient is a model for converting gas to power at a gas power plant.
 In our example, it will be situated between the nodes S$4$ and N$1$ and convert a gas flow $\fuel$ into a real power output $P$ according to
\begin{equation}\label{eq:2}
  \fuel(P) = a_0 + a_1 P + a_2 P^2.
\end{equation}
The flow $\fuel$ must be diverted from the pipeline network, hence the coupling condition at node S$4$ must be changed to
\begin{equation}\label{eq:1}
  \begin{aligned}
    p_\text{in} &= p_\text{out},\\
    q_\text{in} &= q_\text{out}+\fuel \ .
  \end{aligned}
\end{equation}
The details of simulation of such a combined network and the treatment of all arising mathematical issues can be found in~\cite{fokken_goettlich_kolb}.
 We now showcase a concrete example.

\section{Numerical results}
\label{sec:2}

\subsection{Problem setup}
\label{sec:setup}

As already noted, we consider a small part of the GasLib-40 network from [2] consisting of 7 pipelines with a total length of 152km. This network is extended by a compressor station and additionally connected to a power grid with 9 nodes by a gas-to-power generator. For this coupled gas-power network, we simulate a sudden increase in power demand within the power grid and study its effect on the gas network. The considered compressor station is supposed to compensate part of the pressure losses in the gas network such that a given pressure bound is satisfied all the time, while power consumption of the compressor is minimized.

To complete the problem description, the following initial and boundary conditions are given:
  \begin{itemize}
  \item $P(t)$ and $Q(t)$ at load nodes,
  \item $P(t)$ and $V(t)$ at generator nodes,
  \item $V(t)$ and $\phi(t)$ at the slack bus,
  \item $p(t)$ at S$5$,
  \item $q(t)$ at S$25$,
  \item $p(x,0)$, $q(x,0)$ for all pipelines.
  \end{itemize}

\begin{figure}[htb]
  \centering
  \begin{subfigure}{0.49\textwidth}
    \centering
    %\includestandalone[width=\textwidth,height=0.2\textheight]{graphs/demandn5}
    \includegraphics[width=\textwidth,height=0.2\textheight]{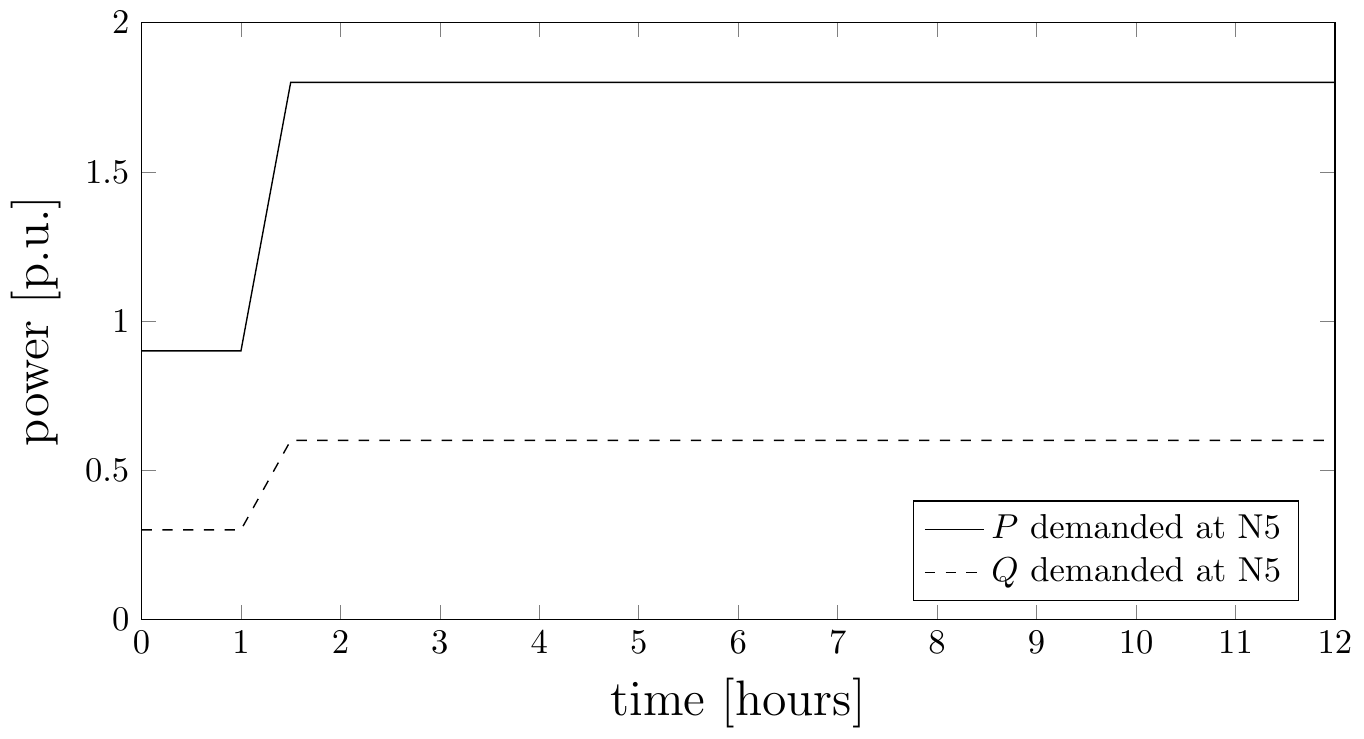}
  \end{subfigure}
  \begin{subfigure}{0.49\textwidth}
    \centering
    %\includestandalone[width=\textwidth,height=0.2\textheight]{graphs/demandslack}
    \includegraphics[width=\textwidth,height=0.2\textheight]{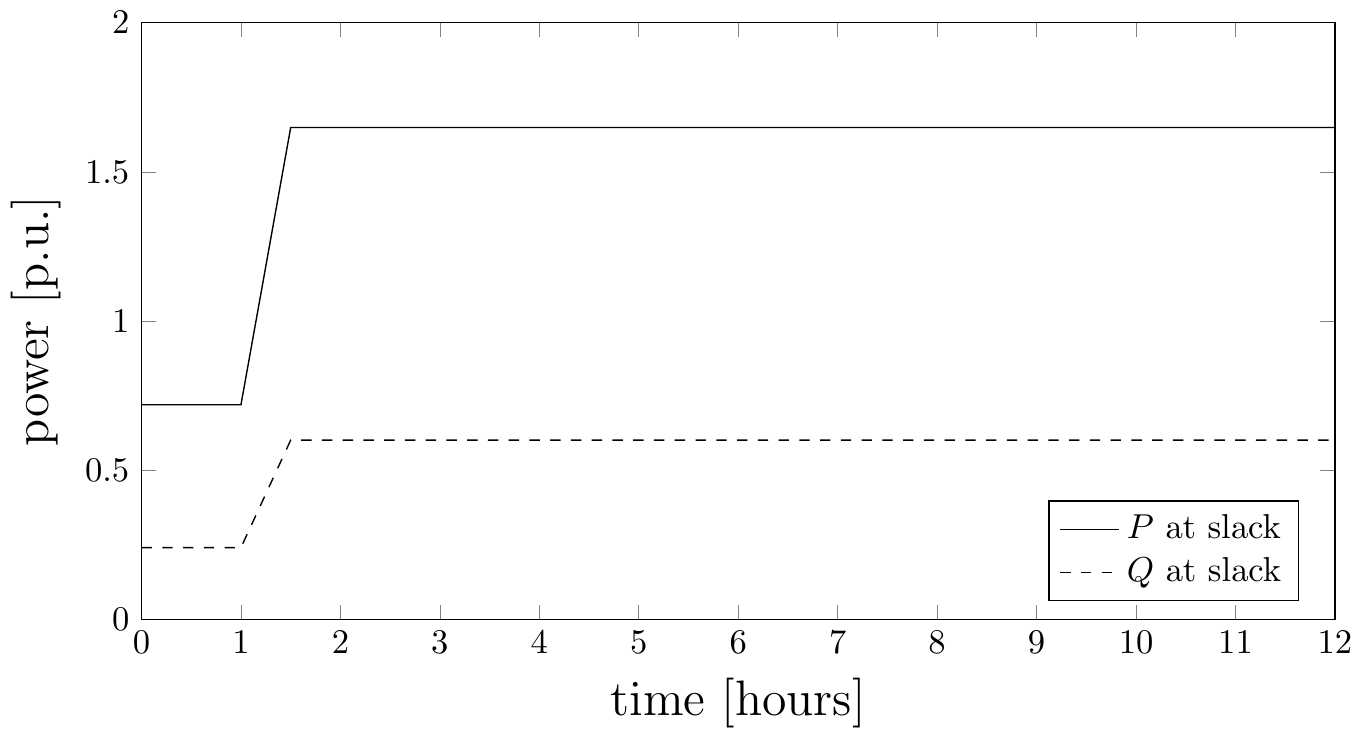}
  \end{subfigure}
  \caption{Power and reactive power at demand node N5 (above) and the slack bus (below).}
  \label{fig:powerN5slack}
\end{figure}

More precisely, the initial conditions for the power network are given in Table~\ref{tab:powerGridInitial}. These remain constant over time except for the power demand at node N$5$, which changes linearly between $t=1$ hour and $t=1.5$ hours from $0.9$ \pu\ to $1.8$ \pu\ for the real power and from $0.3$ \pu\ to $0.6$ \pu\ for the reactive power, see also Figure~\ref{fig:powerN5slack}.

\begin{table}[htb]
 \centering
 \caption{Initial conditions of the power grid (p.u.)}
    \begin{tabular}{@{}lllll}
    \toprule
  Node & $P$ & $Q$ & $\vert V\vert$ & $\va$\\
  \midrule
  N1 & -    & -   & 1 & 0 \\
  N2 & 163  & -   & 1 & - \\
  N3 & 85   & -   & 1 & - \\
  N4 & 0    & 0   & - & - \\
  N5 & -90  & -30 & - & - \\
  N6 & 0    & 0   & - & - \\
  N7 & -100 & -35 & - & - \\
  N8 & 0    & 0   & - & - \\
  N9 & -125 & -50 & - & - \\
  \bottomrule
 \end{tabular}
 \label{tab:powerGridInitial}
\end{table}

For the gas network the incoming pressure at S$5$ is fixed at 60bar, the outflow at S$25$ is fixed at $q=100 \frac{\text{m}^3}{\text{s}}\cdot \frac{\rho_0}{A_e}$, where $\rho_0 = 0.785 \frac{\text{kg}}{\text{m}^3}$ and $A_e = \pi \frac{d_e^2}{4}$. The fuel consumption parameters we use in equation \eqref{eq:2} are given by $a_0=2$, $a_1=5$, $a_2=10$. Since the data for the considered gas and power networks are taken from different sources, the parameters of the gas-to-power generator are chosen in such a way that a significant influence is caused.

Further, the pressure at S$25$ is supposed to satisfy a lower pressure bound of $41$bar, i.e.,
\begin{equation*}
  p_{\text{S}25}(t)\geq 41\text{bar}
\end{equation*}
for all times $t$, where we consider a time horizon of $T = 12$ hours.

As we will see below, the compressor station between nodes S0 and S17 will have to run at a certain time, i.e., $u(t) > 0$, to keep this pressure bound. In general, a solution to the described optimal control problem consists of the control $u(t)$ and the entire network state for all times $t$, i.e.,
  \begin{itemize}
  \item $P(t)$, $Q(t)$, $V(t)$, $\phi(t)$ for all nodes in the power grid,
  \item $p(x,t)$, $q(x,t)$ for all pipelines in the gas network,
  \end{itemize}
fulfilling the given model equations and the pressure constraint.

\subsection{Discretization and Optimization Schemes}
To solve the described optimal control problem, we follow a first-discretize-then-optimize approach. The model equations of the power grid only require a discretization in time, which means that the given boundary conditions and the powerflow equations~\eqref{eq:3} hold for discrete times $t_j = j \Delta t$ with $\Delta t = 15$ minutes and $j \in \{0, \ldots, 48\}$ in our scenario. The discretization of the isentropic Euler equations within the pipelines of the gas network additionally requires a spatial grid (here with grid sizes $\Delta x_e \approx 1$ km) and an appropriate discretization scheme. Here we apply an implicit box scheme \cite{KolbLangBales2010}, which allows the application of large time steps as~$\Delta t = 15$ minutes for the considered spatial grid. Considering the isentropic Euler equations as a system of balance laws of the form
\begin{equation*}
  \label{eq:4}
  y_t + f(y)_x = g(y),
\end{equation*}
the applied scheme reads
\begin{align*}
  \frac{Y^{n+1}_{j-1} + Y^{n+1}_{j}}{2} = &\frac{Y^{n}_{j-1} + Y^{n}_{j}}{2}\\
  &- \frac{\Delta t}{\Delta x}\left( f(Y^{n+1}_j)- f(Y^{n+1}_{j-1}) \right)\\
  &+ \Delta t \frac{g(Y^{n+1}_j)+g(Y^{n+1}_{j-1})}{2}.
\end{align*}
The numerical approximation of the balance law is thought in the following sense:
\begin{align*}
   Y_j^n \approx y(x,t) \quad \text{for} \quad &x\in \big[ (j-\tfrac{1}{2})\Delta x , (j+\tfrac{1}{2})\Delta x \big) ,\\ &t \in \big[ n\Delta t , (n+1)\Delta t \big).
\end{align*}

Together with the algebraic equations modelling the compressor station and the coupling and boundary conditions, the discretization process results in a system of nonlinear equations for all state variables of the coupled gas-power network. For simulation purposes, the entire discretized system is solved with Newton’s method. Note that the system can be solved time-step per time-step and that we exploit the sparsity structure of the underlying Jacobian matrices.

So far we can only compute the state of the considered gas-power network (for a given compressor control $u(t)$) and evaluate quantities of interest like the power consumption of the compressor station or the pressure constraint within the time horizon of the simulation. 
For the given time discretization, the compressor costs~\eqref{eq:compressorCosts} are approximated by the trapezoidal rule and formally contained in $J(u,y(u))$ below.
Next, we want to solve the (discretized) optimal control problem, i.e., to find control values $u(t_j)$ such that the pressure constraint is satisfied, while the power consumption of the compressor is minimized. For this purpose, we apply the interior point optimization code IPOPT~\cite{Ipopt} to the (reduced) discretized optimal control problem. In addition to our simulation procedure to evaluate quantities of interest for a given control, IPOPT further requires gradient information of those quantities with respect to the control. Based on the considered discretization, such information can be efficiently computed by an adjoint approach, which we briefly describe in the following. Therefore, we consider the discretized optimal control problem in the following abstract form:
\begin{equation*}
 \begin{array}{cl}
  \min & J(u,y(u))\\
  \text{s.t.} & u \in \RR^\Nu,\ y \in \RR^\Ny,\\
  & \text{model equations: } E(u,y(u)) = 0,\\
  & \text{constraints: } g(u,y) \ge 0.
 \end{array}
\end{equation*}
The vector $u$ contains all control variables (here the compressor control $u(t_j)$ for all times $t_j = j\Delta t \in [0,24]$) and $y$ contains all state variables of the coupled gas-power network of all time steps $t_j$. The mapping $u \rightarrow y(u)$ is implicitly given by our simulation procedure. Thus, IPOPT does not have to care about the model equations formally summarized in $E(u,y)$, but only about the further constraints $g(u,y)$ and minimizing the objective function $J(u,y)$. Accordingly, we need to provide total derivatives of $J$ and $g$ with respect to $u$.

In the following, we consider the computation of these derivatives only for the objective function, since the procedure is identical for the constraints, and we follow the description given in~\cite{KolbGoettlich2015}. First of all, the chain rule yields
\begin{multline}\label{eq:tdJ}
 \td{u} J(u,y(u)) = \pd{u} J(u,y(u))\\ + \pd{y} J(u,y(u)) \td{u} y(u).
\end{multline}
While the partial derivatives of $J$ with respect to $u$ and $y$ can be directly computed, the derivatives of the states $y$ with respect to the control $u$ are only implicitly given. Differentiating the model equations $E(u,y(u)) = 0$ yields
\begin{equation*}
 0 = \td{u} E(u,y(u)) = \pd{u} E(u,y(u)) + \pd{y} E(u,y(u)) \td{u} y(u)
\end{equation*}
and therewith (formally)
\begin{equation}\label{eq:tdy}
 \td{u} y(u) = - \left(\pd{y} E(u,y(u))\right)^{-1} \pd{u} E(u,y(u)).
\end{equation}
Even though the partial derivatives on the right-hand-side of~\eqref{eq:tdy} can be directly computed, one would have to solve $\Nu$ systems of linear equations here.
Instead of that, we insert~\eqref{eq:tdy} into~\eqref{eq:tdJ} and get
\begin{multline*}
 \td{u} J(u,y(u)) = \pd{u} J(u,y(u))\\ \underbrace{- \pd{y} J(u,y(u)) \left(\pd{y} E(u,y(u))\right)^{-1}}_{= \xi^T} \pd{u} E(u,y(u)).
\end{multline*}
With the so-called adjoint state $\xi$ as the solution of the adjoint equation
\begin{equation}\label{eq:adjointEquation}
  \left(\pd{y} E(u,y(u))\right)^T \xi = - \left(\pd{y} J(u,y(u))\right)^T \,,
\end{equation}
we finally have
\begin{equation}\label{eq:tdJ_final}
 \td{u} J(u,y(u)) = \pd{u} J(u,y(u)) + \xi^T \pd{u} E(u,y(u)) \,.
\end{equation}
It is the fact that~\eqref{eq:adjointEquation} is a single linear system and has a special structure, which can be easily exploited (see for instance~\cite{DissKolb,KolbLang2012}), which makes the computation of derivatives via the presented adjoint approach very efficient.  Nevertheless note that for given control variables $u$ one still has to solve the model equations to get $y(u)$, before one may compute the gradient information via~\eqref{eq:adjointEquation} and~\eqref{eq:tdJ_final}.

\subsection{Results}
\label{sec:results}

\begin{figure}[htb]
  \centering
  \begin{subfigure}{0.49\textwidth}
    \centering
    %\includestandalone[width=\textwidth,height=0.2\textheight]{graphs/inflows5}
    \includegraphics[width=\textwidth,height=0.2\textheight]{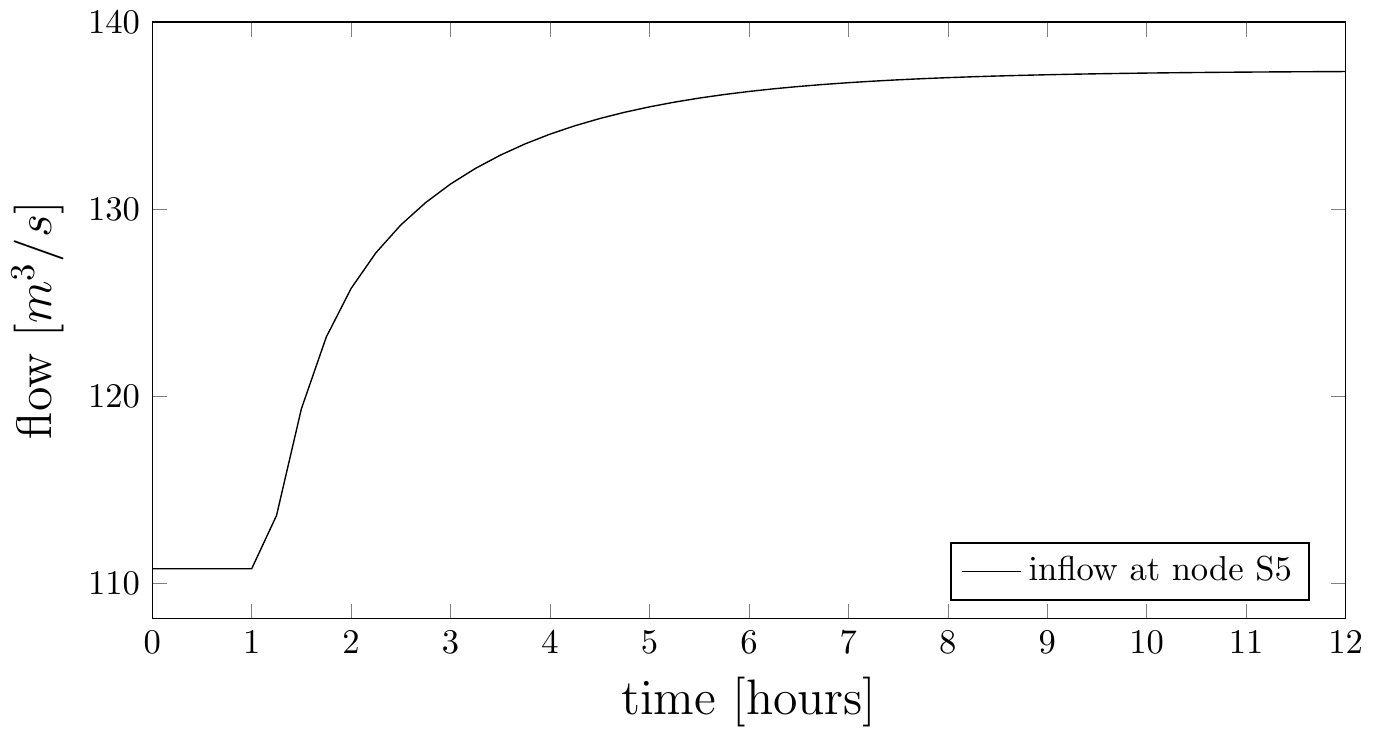}
      \end{subfigure}
    \caption{Inflow at the node S5.}
 \label{fig:inflowS5}
\end{figure}

\begin{figure}[htb]
  \centering
    \begin{subfigure}{0.49\textwidth}
      \centering
        %\includestandalone[width=\textwidth,height=0.2\textheight]{graphs/pressures25}
        \includegraphics[width=\textwidth]{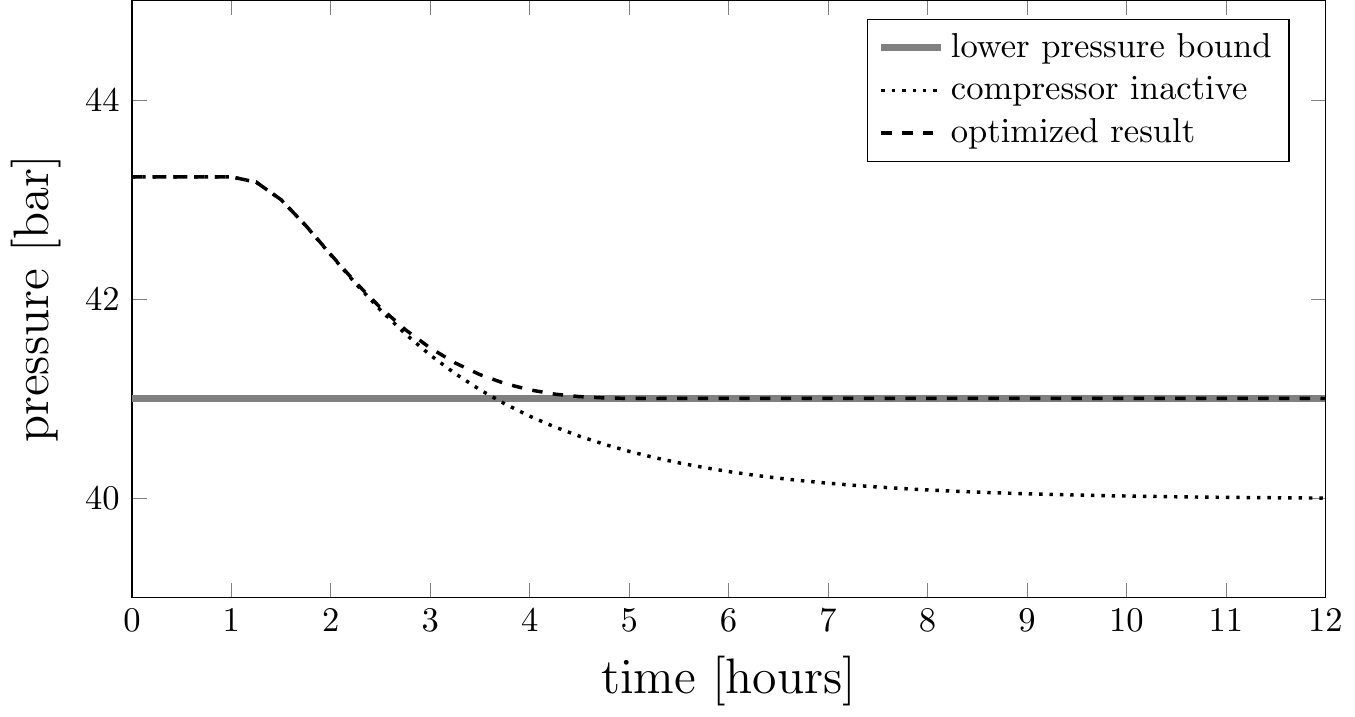}  
      \end{subfigure}
 \caption{Pressure at the node S25.}
 \label{fig:pressureS25}
\end{figure}

We first discuss the simulation with inactive compressor.
 In the course of the simulation, due to the increase in power demand at node N$5$, the power demand at the slack bus rises as well and leads to increased fuel demand at node S$4$.
 This increases the inflow into the gas network, as can be seen in Figure \ref{fig:inflowS5}.
 Also the pressure in the final node S$25$ decreases and violates the pressure bound after approximately 4 hours, see Figure~\ref{fig:pressureS25}.

After the optimization procedure, the compressor station compensates part of the pressure losses in the gas network such that the pressure bound is satisfied all the time. Since the power consumption of the compressor station is minimized within the optimization, the pressure constraint is active after roughly 4 hours (see again Figure~\ref{fig:pressureS25}), i.e., the compressor station applies as little as possible energy.

\section{Conclusion}
The proposed optimization model allows 
to predict pressure transgressions within a coupled gas-to-power framework. 
Simulation and optimization tasks are efficiently solved by exploiting the underlying nonlinear problem structure 
while keeping the full transient regime. 
This makes it possible to track bounds much more accurately than with a steady state model, thereby achieving lower costs.

% \printbibliography[heading=bibintoc]
\bibliographystyle{siam}
\bibliography{bibliography}

%\begin{acknowledgements}
%\end{acknowledgements}


%\bibliography{bibliography}   % name your BibTeX data base

%\begin{thebibliography}{}
% etc
%\end{thebibliography}

\end{document}